\documentclass[11pt,fleqn]{article}
\usepackage{hyperref}
\usepackage{amsfonts}
\usepackage{amssymb}
\usepackage{graphicx}
\newcommand{\ol}{\setlength{\itemsep}{0pt.}\begin{enumerate}}
\newcommand{\eol}{\end{enumerate}\setlength{\itemsep}{-\parsep}}
\newcommand{\ignore}[1]{}
\newtheorem{THEOREM}{Theorem}[section]
\newenvironment{theorem}{\begin{THEOREM} \hspace{-.85em} {\bf :}
}%
                        {\end{THEOREM}}
\newtheorem{LEMMA}[THEOREM]{Lemma}
\newenvironment{lemma}{\begin{LEMMA} \hspace{-.85em} {\bf :} }%
                      {\end{LEMMA}}
\newtheorem{COROLLARY}[THEOREM]{Corollary}
\newenvironment{corollary}{\begin{COROLLARY} \hspace{-.85em} {\bf
:} }%
                          {\end{COROLLARY}}
\newtheorem{PROPOSITION}[THEOREM]{Proposition}
\newenvironment{proposition}{\begin{PROPOSITION} \hspace{-.85em}
{\bf :} }%
                            {\end{PROPOSITION}}
\newtheorem{DEFINITION}[THEOREM]{Definition}
\newenvironment{definition}{\begin{DEFINITION} \hspace{-.85em} {\bf
:} \rm}%
                            {\end{DEFINITION}}
\newtheorem{EXAMPLE}[THEOREM]{Example}
\newenvironment{example}{\begin{EXAMPLE} \hspace{-.85em} {\bf :}
\rm}%
                            {\end{EXAMPLE}}
\newtheorem{CONJECTURE}[THEOREM]{Conjecture}
\newenvironment{conjecture}{\begin{CONJECTURE} \hspace{-.85em}
{\bf :} \rm}%
                            {\end{CONJECTURE}}
\newtheorem{MAINCONJECTURE}[THEOREM]{Main Conjecture}
\newenvironment{mainconjecture}{\begin{MAINCONJECTURE} \hspace{-.85em}
{\bf :} \rm}%
                            {\end{MAINCONJECTURE}}
\newtheorem{PROBLEM}[THEOREM]{Problem}
\newenvironment{problem}{\begin{PROBLEM} \hspace{-.85em} {\bf :}
\rm}%
                            {\end{PROBLEM}}
\newtheorem{QUESTION}[THEOREM]{Question}
\newenvironment{question}{\begin{QUESTION} \hspace{-.85em} {\bf :}
\rm}%
                            {\end{QUESTION}}
\newtheorem{REMARK}[THEOREM]{Remark}
\newenvironment{remark}{\begin{REMARK} \hspace{-.85em} {\bf :}
\rm}%
                            {\end{REMARK}}
\newtheorem{CLAIM}[THEOREM]{Claim}
\newenvironment{claim}{\begin{CLAIM} \hspace{-.85em} {\bf :}
\rm}%
                            {\end{CLAIM}}
\newcommand{\thm}{\begin{theorem}}
\newcommand{\lem}{\begin{lemma}}
\newcommand{\pro}{\begin{proposition}}
\newcommand{\dfn}{\begin{definition}}
\newcommand{\rem}{\begin{remark}}
\newcommand{\xam}{\begin{example}}
\newcommand{\cnj}{\begin{conjecture}}
\newcommand{\mcnj}{\begin{mainconjecture}}
\newcommand{\prb}{\begin{problem}}
\newcommand{\que}{\begin{question}}
\newcommand{\cor}{\begin{corollary}}
\newcommand{\clm}{\begin{claim}}
\newcommand{\prf}{\noindent{\bf Proof:} }
\newcommand{\ethm}{\end{theorem}}
\newcommand{\elem}{\end{lemma}}
\newcommand{\epro}{\end{proposition}}
\newcommand{\edfn}{\end{definition}}
\newcommand{\erem}{\bbox\end{remark}}
\newcommand{\exam}{\bbox\end{example}}
\newcommand{\ecnj}{\bbox\end{conjecture}}
\newcommand{\emcnj}{\bbox\end{mainconjecture}}
\newcommand{\eprb}{\bbox\end{problem}}
\newcommand{\eque}{\bbox\end{question}}
\newcommand{\ecor}{\end{corollary}}
\newcommand{\eclm}{\end{claim}}

\newcommand{\beqn}{\begin{equation}}
\newcommand{\eeqn}{\end{equation}}

\newcommand{\bbox}{\begin{flushright} $\Box $ \end{flushright}}
\newcommand{\qed}{\bbox}

\def\eps{\epsilon}

\def \1{\mathbf 1}
\newcommand{\X}[0]{{\cal X}}
 \newcommand{\Y}[0]{{\cal Y}}

\def\<{\left<}
\def\>{\right>}
\def \({\left(}
\def \){\right)}

\def \8{\infty}

\newcommand{\remove}[1]{}

\newcommand{\enote}[1]{{\bf (Ehud:} {#1}{\bf ) }}

\title{An information-theoretic proof of a hypercontractive inequality}
\author{Ehud Friedgut\footnote{Weizmann Institute of Science, Israel. Research supported in part by I.S.F. grant 0398246, BSF grant 2010247 and MINERVA grant 712023.}}
\begin{document}
\maketitle
\begin{abstract}
The famous hypercontractive estimate discovered independently by Gross \cite{Gross}, Bonami \cite{Bonami} and Beckner \cite{Beckner}, has had great impact on combinatorics and theoretical computer science since first used in this setting in  
the seminal KKL paper \cite{KKL}. The usual proofs of this inequality begin with the two-point space where some elementary calculus is used, and then generalise immediately by induction on the dimension, using submultiplicativity (Minkowski's integral inequality). In this paper we prove the inequality using information theory. We compare the entropy of a pair of correlated vectors in $\{0,1\}^n$ to their separate entropies, analysing them bit by bit (not as a figure of speech, but as the bits are revealed), using the chain rule of entropy.
\end{abstract}

\section{Introduction}
The inequality that we consider in this note is a two-function version of a famous hypercontractive inequality due, independently, to Gross \cite{Gross}, Bonami \cite{Bonami} and Beckner \cite{Beckner}. This inequality, first introduced to the combinatorial landscape in the seminal KKL paper \cite{KKL}, has become one of the cornerstones of the analytical approach to Boolean functions and theoretical computer science, see e.g.  \cite{BKS}, \cite{DFR}, \cite{FKN}, \cite{KO},\cite{MORSS}, \cite{MosselArrow}, and many many others. See chapter 16 of \cite{O} for a historical background. 

Let $\eps \in (0,1)$. We will be considering an operator $T_\eps$ which acts on real valued functions  on $\{0,1\}^n$. We consider two equivalent definitions of the operator; A spectral definition and a more probabilistic one. The first definition is via the eigenfunctions and eigenvalues of the operator. The eigenfunctions are precisely the Walsh-Fourier characters,
 $\{ u_X \}_{X \in \{0,1,\}^n}$, which form a complete orthonormal system under the standard inner product on $\{0,1\}^n$. We recall the definition of these characters. For $X,Y \in \{0,1\}^n$
$$
u_X(Y)=(-1)^{\sum X_i Y_i}.
$$
Given a function $f : \{0,1\}^n \rightarrow \mathbb{R}$, and its unique Fourier expansion, $f = \sum \hat{f}(X) u_X$ the action of $T_{\epsilon}$ on $f$ is defined by
$$
T_\eps(f) = \sum \eps^{\sum X_i}\hat{f}(X) u_X.
$$
This definition of $T_\eps$ stresses the fact that it "focuses'' on the low-frequency part of the Fourier spectrum, an idea that was a crucial element in the KKL proof \cite{KKL}. 

For the other definition let $X,Y$ be random variables taking values in $\{0,1\}^n$. Either for fixed $X$, or any distribution of $X$, let $Y$ be such that for every coordinate $i$, independently, $Y_i$ is chosen so that $Pr[X_i=Y_i]= \frac{1+\eps}{2}$. (Or, if one prefers the $\{-1,1\}^n$ setting, the restriction is $E[X_iY_i] = \eps$.) Note that if $X$ is chosen uniformly, then the marginal distribution of $Y$ is also uniform. We call such a pair 
$(X,Y)$ an  $\epsilon${\em-correlated pair}. Then one can define
$$
T_{\eps}(f)(X) = E[f(Y)|X].
$$
It is not hard to verify that these two definitions of $T_\eps$ are equivalent.  
The second definition, which is the one we will be working with in this paper, stresses the connection of this operator to random walks and isoperimetric inequalities, as it enables one to bound the probability of an $\eps$-correlated pair of random points $X,Y$ to belong to given sets. It also explains the fact (that will be made formal shortly) that $T_\eps(f)$ is smoother than $f$, as the operator is an averaging operator.
Without further ado, here is the statement of the inequality.
\thm[Gross,Bonami,Beckner]
Let $f : \{0,1\}^n \rightarrow  \mathbb{R}^{\ge 0}$, and let $\epsilon \in [0,1]$. Then
\beqn\label{88}
|T_\eps(f)|_2 \leq |f|_{1+\eps^2},
\eeqn
where, as usual, $|g|_p = E[|g|^p]$, and the expectation is with respect to the uniform measure.
\ethm
{\bf Remarks:} 
\begin{itemize}
\item There are various refinements of this inequality either dealing with non-uniform measure, the norm of $T_\eps$ as an operator from $L^p$ to $ L^q$ for $q \not = 2$, studying products of a base space with more than two points, and also
a reverse inequality that deals with the case $p,q < 1$. See \cite{Olesz},\cite{Wolf},\cite{MORSS}, \cite{Borel}. It would not be surprising if the method of this note could be extended to cover such cases too.
\item It's not difficult to see that (\ref{88}) is equivalent to the following. 
Let $f, g :  \{0,1\}^n \rightarrow \mathbb{R}$, let $X$ be uniformly distributed on $\{0,1\}^n$, and let $X,Y$ be an $\epsilon$-correlated pair. Then 
\beqn\label{77}
E[f(X)g(Y)] \leq |f|_{1+\eps}|g|_{1+\eps}.
\eeqn
This is the inequality proven in this paper.
 \item A major portion of the applications of the hypercontractive inequality deal with the case when $f$ and $g$ are Boolean functions. We will start our proof with this setting, and then show how a small variation deals with the general case.
\end{itemize}
In this paper we apply an information-theoretic approach to proving (\ref{77}) for Boolean functions, trying to analyse the pair $(X,Y)$ as the coordinates of $X$ and $Y$ are revealed to us one by one. Since all known direct proofs of (\ref{88}) use induction, it is not surprising that one should adopt such a sequential approach. The difference is that the usual proofs begin with the two point space, and proceed by induction, using  submultiplicativity of the product operator and Minkowski's integral inequality, whereas we use the chain rule of entropy, exposing the bits of the vectors in question one by one, and comparing the amount of information of their joint distribution with the information captured by their marginal distributions. Fortunately it turns out that regardless of the prefixes revealed so far, at every step the conditional entropies obey the same inequality.  

This is not the first application of entropy to this hypercontractive inequality. 
The connection between hypercontractivity-type inequalities and information theory is perhaps first addressed in \cite{AG}.
In \cite{FR} the dual form of the hypercontractive inequality is proven for the case of comparing the 2-norm and the 4-norm of a low degree polynomial on $\{0,1\}^n$. Blais and Tan, \cite{BT}, managed to improve this approach and, surprisingly,  extract the precise optimal hypercontractive constant for comparing the 2-norm and the $q$-norm of such polynomials, for all positive even integers $q$. Both these proofs analyse the Fourier space rather than the primal space - and use no induction at all.

One final remark regarding the proof in this paper. Although it is not difficult, it is probably, to date, the most involved proof of the inequality from a technical point of view. I believe that nonetheless it is worthwhile to add it to the list of existing proofs, because it offers a new point of view which directly addresses the notion of studying projections of the joint distribution of  a pair of $\epsilon$-correlated vectors. 

\section{Main Theorem}
\subsection{The Boolean case}
\thm
Let $\eps \in (0,1)$, and let $\X, \Y \subseteq \{0,1\}^n$ be nonempty. Let $X$ be uniformly distributed on $\{0,1\}^n$, and let $Y$ be such that for each $1 \le i \le n$ independently $Pr[X_i = Y_i] = \frac{1+\eps}{2}$.
Then
\beqn\label{1}
E[{\bf 1}_{\X}(X){\bf 1}_{\Y}(Y)] \le \left(\mu(\X) \mu(\Y)\right) ^{\frac{1}{1+\eps}} .
\eeqn
 \ethm
 \remove{Remark: When we take $\X=\Y$ and $f={\bf 1}_{\X}$ this gives (\ref{1}) in the form 
 $$
 \sum\hat{f}^2(S)\eps^{|S|} \leq |f|_{1+\epsilon}^2,
 $$
 (which looks more familiar  with $\eps^2$ instead of $\eps$). 
 }
   \prf
 For $X,Y \in \{0,1\}^n$ let $a(X,Y)$ denote the number of coordinates on which $X$ and $Y$ agree, and $d(X,Y)$ be the number of coordinates on which they differ.
 Then the theorem is equivalent (by straightforward manipulation) to
 \beqn\label{2}
 \log \left( \sum_{X \in \X} \sum_{Y \in \Y}  (1+\eps)^{a(X,Y)}(1-\eps)^{d(X,Y)} \right) \leq
  \frac{1}{1+\eps} \left(  2\eps n + \log (|\X|) + \log(|\Y|)       \right),
 \eeqn
 where all logs are base 2. As usual in proofs using entropy, it suffices, by continuity, to treat the case where $\eps$ is rational. Let $s \leq r$ be natural numbers such that $\frac{1+\eps}{2}=\frac{r}{r+s}, \frac{1-\eps}{2}=\frac{s}{r+s}$.
 Then (2) reduces to
 \beqn\label{3}
 \log \left( \sum_{X \in \X} \sum_{Y \in \Y}  r^{a(X,Y)}s^{d(X,Y)} \right) \leq
  n(\log(r+s) -s/r) + \frac{r+s}{2r} (\log (|\X|) + \log(|\Y|)).       
 \eeqn
We will express the left hand side of this expression as the entropy of a random variable, and proceed to expand it according to the chain rule. First let $A_{00}, A_{11} , A_{10}, A_{01}$ be four disjoint sets with
$$
|A_{00}|=|A_{11}| = r, |A_{01}|=|A_{10}| = s.
$$  
Next let $(X,Y,Z)$ be a random variable which is distributed uniformly over all triples such that $X \in \X, Y \in \Y,$ and for every $1 \leq i \leq n$, $Z_i \in A_{X_iY_i}.$ Clearly $Z$ determines $X$ and $Y$ so that 
$$
H(X,Y,Z) = H(Z) = 
 \log \left( \sum_{X \in \X} \sum_{Y \in \Y}  r^{a(X,Y)}s^{d(X,Y)} \right).
 $$
 Next, for any vector $W \in \{0,1\}^n $ we denote $(W_1,\ldots, W_{i-1}) := W_{<i}$. So by the chain rule we have
 $$
 H(Z) = \sum_i H(Z_i | Z_{<i})
 $$
 and 
 $$
 H(X)= \sum_i H(X_i | X_{<i}) \leq \log (|\X|),  H(Y)= \sum_i H(Y_i | Y_{<i}) \leq \log (|\Y|).
 $$
 Hence it suffices to prove 
 \begin{equation}\label{eq:entropy}
  H(Z ) \leq  
   \frac{r+s}{2r}\left( H(X)
  +  H(Y) \right)
   +n( \log(r+s)  - s/r)  .
\end{equation}

 Noting that any fixed value of $Z_{<i}$  determines the values  $X_{<i}$ and  $Y_{<i}$ we have
 $$
H(X_i | Z_{<i}) \leq H(X_i | X_{<i}),
 $$
 and
 $$
H(Y_i | Z_{<i}) \leq H(Y_i | Y_{<i}),
 $$
so inequality (\ref{eq:entropy}) will follow from  
 \begin{claim}
 Denote a fixed $Z_{<i} := Past$. Then for all fixed values of $Past$
  $$
  H(Z_i | Past) \leq  
   \frac{r+s}{2r}\left( H(X_i |  Past)
  +  H(Y_i |Past )  \right)
   + \log(r+s)  - \frac{s}{r}  
 $$
 \end{claim}
 \prf
We condition on $Past$, and by abuse of notation drop the dependency in the notation,
e.g. $H(X)$ and $ H(Z|X)$ rather than $H(X| Past) , H(Z| X, Past)$, etc. We also drop the index $_i$ and write $X$ for $X_i$ etc., so, using the new notation, we want to prove (for all integers $ r \ge s \ge 0$)
  $$
  H(Z ) \leq  
   \frac{r+s}{2r}\left( H(X)
  +  H(Y) \right)
   + \log(r+s)  - \frac{s}{r}  
 $$

Note that $H(Z) = H(X,Y) + H(Z|X,Y) = H(X,Y) + \log r (Pr [X=Y]) + \log s (Pr[X \not = Y])$, so we need to prove

$$
 \frac{r+s}{2r}\left( H(X)
  +  H(Y) \right) - H(X,Y) - \log r (Pr [X=Y]) - \log s (Pr[X \not = Y])
   + \log(r+s)  - \frac{s}{r} \ge 0
$$ 
Since this expression is invariant when $r$ and $s$ are multiplied by any positive constant we can set $r=1$ and denote $\delta := s/r$. Next, for a joint distribution of $X$ and $Y$ on $\{0,1\}^2$, and $(i,j) \in \{0,1\}^n$ let 
$P_{ij}= Pr[X=i, Y=j]$  
So we want to prove
\begin{equation}\label{eq:defineF}
F_\delta\left(
\begin{array}{cc}
P_{01}  & P_{11}    \\
P_{00}  & P_{10}     
\end{array}
\right)
 :=   \frac{1+\delta}{2}\left( H(X)
  +  H(Y) \right) - H(X,Y)  - \log \delta (Pr[X \not = Y])
   + \log(1+\delta)  - \delta \ge 0
\end{equation}
We know (and can check directly from the formula) that equality holds when $\X=\Y=\{0,1\}^n$ in which case we have
\begin{equation}\label{eq:solution}
\left(
\begin{array}{cc}
P_{01}  & P_{11}    \\
P_{00}  & P_{10}     
\end{array}
\right)
=
\left(
\begin{array}{cc}
\frac{\delta}{2+2\delta}  & \frac{1}{2+2\delta}    \\
\frac{1}{2+2\delta} & \frac{\delta}{2+2\delta}  
\end{array}
\right).
\end{equation}
We wish to show that this is the unique minimum. 
To simplify notation (and save indices) we denote 
\[
\left(
\begin{array}{cc}
P_{01}  & P_{11}    \\
P_{00}  & P_{10}     
\end{array}
\right)
:=
\left(
\begin{array}{cc}
a  & b    \\
c & d  
\end{array}
\right),
\]
and attempt to minimise 
$F_\delta\left(
\begin{array}{cc}
a  & b    \\
c  & d     
\end{array}
\right)
$
under the constraints 
 $$
 a,b,c,d \ge 0 \mbox{ and } a+b+c+d=1
 .$$
Using Lagrange multipliers we deduce (after some simple cancelations) that at a local minimum in the interior of the region in question one must have 
\begin{equation}\label{eq:a}
\frac{\delta}{a}\left[{(a+b)(a+c)}\right]^{(1+\delta)/2} = 
\end{equation}
\begin{equation}\label{eq:b}
\frac{1}{b}\left[{(a+b)(b+d)}\right]^{(1+\delta)/2}=
\end{equation}
\begin{equation}\label{eq:c}
\frac{1}{c}\left[(a+c)(c+d)\right]^{(1+\delta)/2}=
\end{equation}
\begin{equation}\label{eq:d}
\frac{\delta}{d}\left[(c+d)(b+d)\right]^{(1+\delta)/2}. 
\end{equation}
From the fact that (\ref{eq:a})$\cdot$(\ref{eq:d})=(\ref{eq:b})$\cdot$(\ref{eq:c}) we get that 
\begin{equation}\label{eq:ratio}
ad = \delta^2bc
\end{equation}
Next we plug (\ref{eq:ratio}) and the restriction $a+b+c+d=1$ into the equation (\ref{eq:a})=(\ref{eq:d}). This yields
\begin{equation}\label{eq:b=c}
\left(\frac{d}{a}\right)^{(1-\delta)/2} = \left[\frac{ 1 +(\delta^{-2}-1)a}{1+(\delta^{-2}-1)d}     \right]^{(1+\delta)/2}.
\end{equation}
Note that for every fixed value of $b+c$, the value of $a+d$ is fixed, so letting $d$ grow from 0 to $1-b-c$, as $a=1-b-c-d$ decreases from $1-b-c$ to 0, we see that 
 the left hand side of (\ref{eq:b=c}) is increasing and the right hand side is decreasing, hence there exists a single solution, which, by inspection, is $a=d$. 
 
 We now know that $a=d$ and that $bc = \delta^{-2} a^2$ and $b+c = 1-2a$. This gives $b$ and $c$ as the roots of the quadratic equation $X^2 - (1-2a)X + \delta^{-2} a^2=0$. Plugging these roots into the equation  (\ref{eq:b})=(\ref{eq:c}), and using $a=d$ yields the following equation for $a$:
 \begin{equation}\label{eq:finala}
 \frac{[1-2a+S(a)][1-S(a)]^{1+\delta}}{[1-2a-S(a)][1+S(a)]^{1+\delta}} =1,
 \end{equation}
 where $S(a)$ denotes $\sqrt{1-4a+4(1-\delta^{-2})a^2}$.
 Now, $a$ can take on values between 0 and $1/2$, as long as $S(a) \ge 0$, so the relevant range is $0\le a \le \frac{\delta}{2+2\delta}$.
 When $a= \frac{\delta}{2+2\delta}$, as required, then $S(a) =0$ and equation (\ref{eq:finala}) clearly holds.
 An elementary calculation shows that for all $\delta \in [0,1]$ the left hand side of (\ref{eq:finala}) is a decreasing function of $a$ in the interval $[0, \frac{\delta}{2+2\delta}]$, so $(a,b,c,d)$ is determined, and there is a unique internal minimum in the region which we are exploring. (The meticulous reader may check that the derivative of the  left hand side of (\ref{eq:finala}) according to $a$ is 
 
 $$
 \frac{16 a^2 (1-\delta) \left(1-S(a)\right)^\delta \left(\1+S(a)\right)^{-\delta-2} (2 a (\delta+1)-\delta)}{\delta^3 \left(S(a)+2 a-1\right)^2 S(a)},
 $$
 whose sign is determined by  $ 2 a (\delta+1)-\delta$, which is negative for all $a \in  (0, \frac{\delta}{2+2\delta})$. )
 
\remove{ *** What's the range of delta? ****}
 
What about points on the boundary of the region? We claim that there can be no minima with negative values on the boundary.
We deal with two cases. The case of 
$\left(
\begin{array}{cc}
a  & b    \\
c & 0   
\end{array}
\right)
$ 
with $a>0$ (and the different rotations of this), and the case of
$\left(
\begin{array}{cc}
0  & b    \\
c & 0  
\end{array}
\right)
$ and its rotation (note that 
$
F_\delta\left(
\begin{array}{cc}
0  & y    \\
z & 0  
\end{array}
\right) \le
$
$F_\delta\left(
\begin{array}{cc}
y  & 0    \\
0 & z    
\end{array}
\right) 
$

In the first case, there is a nearby point in the interior of the region, where the value of $F_\delta$ is smaller, this is because the derivative of
$$
F_\delta\left(
\begin{array}{cc}
a-t  & b    \\
c & t     
\end{array}
\right)
$$ with respect to $t$ at $t=0$ is minus infinity (as one of the summands being derived is $t \log(t)$.)
In the second case the value of $F_\delta$ is non negative:
$$
F_\delta\left(
\begin{array}{cc}
0 & b    \\
1-b & 0     
\end{array}
\right) =
   \delta(b\log(1/b)+(1-b)\log(1/(1-b)))
   + \log(1+\delta)  - \delta, 
$$
which is non-negative, as $\log_2(1+\delta)\ge\delta$ for all $\delta \in [0,1]$.

\remove{
Let us study the faces of the polytope we're interested in by decreasing dimension, i.e. does the matrix 
$\left(
\begin{array}{cc}
a  & b    \\
c & d     
\end{array}
\right)
$ 
have one zero entry, two zeroes, or three.
\\ {\bf One zero.} Noting that 

$F_\delta\left(
\begin{array}{cc}
x  & y    \\
z & 0     
\end{array}
\right) =
$
$F_\delta\left(
\begin{array}{cc}
0  & y    \\
z & x     
\end{array}
\right) \le
$
$F_\delta\left(
\begin{array}{cc}
y  & 0    \\
z & x    
\end{array}
\right) =
$
$F_\delta\left(
\begin{array}{cc}
y  & z    \\
0 & x    
\end{array}
\right)
$
we can restrict to the case
$F_\delta\left(
\begin{array}{cc}
a  & b   \\
c & 0     
\end{array}
\right).
$
\\  Using Lagrange multipliers, and the identity $a+b+c=1$ we get that  necessary condition for an internal minimum
is
$$
(1-b)^{\frac{1+\delta}{2}}b^{\frac{1-\delta}{2}}=(1-c)^{\frac{1+\delta}{2}}c^{\frac{1-\delta}{2}}.
$$
Note that the expression on the LHS is a geometric average of $b$ and $1-b$

can have either no zeros (an interior point, the case we have already covered), two zeros in the same row or column, or three zeroes, i.e it cannot have a single zero, or zeroes only on one of the diagonals. The reason is that these probabilities signify correlations between sets of vectors, and 
can be equal to zero only if one of the sets is empty.
The case of three zeroes boils down to noting that $\log_2(1+\delta) - \delta$ is non-negative for $\delta \in [0,1]$.
Let's consider the case of two zeroes, say
 $\left(
\begin{array}{cc}
P_{01}  & P_{11}    \\
0  & 0
\end{array}
\right).
$
In this case one can check that an infinitesimal change 
$\left(
\begin{array}{cc}
-\Delta  & 0    \\
0  & \Delta     
\end{array}
\right)
$
induces a change of $\frac{\delta -1 +o(1)}{2}\Delta \log(1/\Delta)$ in the value of $F_\delta$ (where the $o(1)$ notation is as $\Delta$ goes to zero). Since this is negative there cannot be a local minimum of this form. 
}
\subsection{The general (non-Boolean) case}
The general case is actually a minor extension of the Boolean one, that follows by adding one more coordinate to each of the random variables $X$ and $Y$.
\thm\label{nBmain}
Let $\eps \in (0,1)$, and let $f,g : \{0,1\}^n \rightarrow \mathbf{R}^{\ge 0}$ . Let $X$ be uniformly distributed on $\{0,1\}^n$, and let $Y$ be such that for each $1 \le i \le n$ independently $Pr[X_i = Y_i] = \frac{1+\eps}{2}$.
Then
\beqn\label{nB}
E[f(X)g(Y)] \le |f|_{1+\epsilon} |g|_{1+\epsilon}.
\eeqn
 \ethm 
 \prf
By continuity it suffices to consider the case where all values of $f$ and $g$ are rational, and by homogeneity, we can clear common denominators and assume that they are integer valued. Now we wish to prove a slight extension of (\ref{3}), namely
 \begin{equation}\label{LHS}
\log \left( \sum_{X \in \X} \sum_{Y \in \Y}  r^{a(X,Y)}s^{d(X,Y)}f(X)g(Y) \right) \leq
\end{equation}
\begin{equation}\label{RHS}
 n(\log(r+s) -s/r) + \frac{r+s}{2r} \left( \log \left(\sum_{X \in {0,1}^n} f(X)^{\frac{2r}{r+s}}\right) +  \log \left(\sum_{Y \in {0,1}^n} f(Y)^{\frac{2r}{r+s}}\right) \right). 
\end{equation}

To this end we now define as before the random variables $X,Y,Z$ and add two more integer random variables $a$ and $b$,
and take $(X,Y,Z,a,b)$ uniformly, with  $(X,Y,Z)$ as before and the additional constraint that $a \in \{1,\ldots, f(X)\}, b \in \{1,\ldots, g(Y)\}$.
Now (\ref{LHS}) is precisely $H(Z,a,b) = H(Z) + H(a|Z) + H(b|Z) = H(Z) + E_X[\log(f(X))] + E_Y[\log(g(Y))]$. On the other hand, note that 
for any function $t$ it holds that 
$$
H(X) + E_X[\log(t(X))] \le \log\left( \sum_X t(X) \right) 
$$
In particular 
$$
\frac{r+s}{2r}(H(X))+ E_X[\log(f(X))] \le \frac{r+s}{2r} \left(\log \left(\sum_{X \in {0,1}^n} f(X)^{\frac{2r}{r+s}}\right)\right)
$$
and
$$
\frac{r+s}{2r}(H(Y))+ E_Y[\log(g(Y))] \le \frac{r+s}{2r} \left(\log \left(\sum_{Y \in {0,1}^n} g(Y)^{\frac{2r}{r+s}}\right)\right)
$$
To complete the proof of theorem \ref{nBmain} we just add $ E_X[\log(f(X))] +   E_Y[\log(g(Y))] $ to both sides of the main inequality that we proved in the Boolean case, namely
$$
  H(Z ) \leq  
   \frac{r+s}{2r}\left( H(X)
  +  H(Y) \right)
   +n( \log(r+s)  - s/r)  
$$ 
and we're done.
\qed  
\remove{ \enote{ For skew $(p,q)$-product the joint distribution for $x,y$ is $p^2+\eps p q , (1-\eps)pq, (1-\eps)pq, q^2+ \eps p q$.}
 }

{\bf Acknowledgments:}
I would like to thank David Ellis and Gideon Schechtman for useful conversations, and the anonymous referee for pointing out some inaccuracies, which, once fixed, led to a simplification of the proof. 
I would also like to thank Chandra Nair for spotting some errors in an early draft of this paper, and for alerting me to the existence of two papers, \cite{Nair} and \cite{CCE}, which also adopt an information theoretic approach to hypercontractivity, and also of \cite{NW1} and \cite{NW2} which address similar problems (and refer to an earlier arxiv-version of this paper.)

 \end{document}